\documentclass[11pt]{amsart}
\usepackage[shortalphabetic]{amsrefs}
\usepackage[utf8]{inputenc}
\usepackage{amsthm}
\usepackage[utf8]{inputenc}
\usepackage{appendix}
\usepackage{cite}
\usepackage{amssymb,amsmath,amsthm}
\usepackage{hyperref}

\newcommand{\amsprimary}[1]{{\footnotesize\noindent AMS 2010 \textit{Mathematics subject
classification:} Primary #1\vspace{1pc}}}
\newcommand{\keywordsnames}[1]{{\footnotesize\noindent\textit{Key words:} #1\vspace{1pc}}}

\newtheorem{theorem}{Theorem}

\theoremstyle{definition}

\theoremstyle{remark}

\title[]{On the behaviour of harmonic functions on Riemannian cones}
\author{Jean C. Cortissoz }
\email{jcortiss@uniandes.edu.co}
\address{Department of Mathematics, Universidad de los Andes, Bogot\'a DC, Colombia}
\date{}

\begin{document}

\maketitle

\begin{abstract}
We discuss the behavior of harmonic functions on Riemannian cones as defined below and
Lioville's theorem.

\end{abstract}

\keywordsnames{Liouville's theorem; bounded harmonic functions}

{\amsprimary {31C05, 53C21}}

\tableofcontents

\section{Introduction}

Liouville's theorem states that a bounded entire function 
(that is, a holomorphic function with domain $\mathbb{C}$) must be constant. Furthermore,
any entire function of sublinear growth must be constant: this statement is sharp, as the
function $f\left(z\right)=z$ is of linear growth. since the real and imaginary parts
of a holomorphic function are harmonic, Liouville's theorem in the plane can be stated
as follows

\medskip
Let $u:\mathbb{R}^2\longrightarrow \mathbb{R}$ be a harmonic function. 
If $u$ is bounded, then it is constant (classical statement of Liouville's theorem).
Furthermore, if
$u=o\left(\left|x\right|\right)$ then $u$ is constant (sublinear version
of Liouville's theorem).

\medskip
In fact, Liouville's theorem as stated above, in both versions, is valid in any dimension. There are
several proofs of this fact, and a beautiful proof due to E. Nelson for the classical statement
of Liouville's theorem is given in 
\cite{Nelson}. It is not difficult to extend
Nelson's proof to show the
sublinear version of Liouville's Theorem.

\medskip
Now, if $\left(M,g\right)$ is a Riemannian manifold, the Laplace operator can be defined
in a natural way, and in fact can be written, using 
local coordinates (and assuming Einstein's summation convention) as
\[
\Delta_g f= \frac{1}{\sqrt{\mbox{det}\left(g\right)}}\partial_i\left(\sqrt{\mbox{det}\left(g\right)}g^{ij}\partial_j f\right).
\]
Thus the space of harmonic functions with domain $M$ can be defined in the obvious way: $u\in C^2\left(M\right)$ 
is harmonic if $\Delta_g u=0$. 

The
study of the behaviour of harmonic functions defined on a Riemannian manifold
has been an important subject of study in Geometric Analysis, and in this
study, Liouville's theorem occupies a special place. It was Yau in \cite{Y}
who started it all by extending Liouville's theorem to manifolds of
nonnegative Ricci curvature:
\begin{theorem}
    Let $\left(M,g\right)$ be a Riemannian manifold with $Ric\left(g\right)\geq 0$.
    Let $O$ be a fixed point on $M$ and for $x\in M$ define $r_x=d\left(O,x\right)$,
    the Riemannian distance from $x$ to $O$.
    Let $u:M\longrightarrow \mathbb{R}$ be a harmonic function. If $u=o\left(r_x\right)$, then
    $u$ is constant.
\end{theorem}
Yau's proof is via a gradient estimate, which is along the spirit
of the proof given in an elementary class in Complex Variables. For nonnegative curvature,
Cheng and Yau's result is quite satisfactory, but it leaves open the question
on what happens in the presence of negative curvature.

\medskip
Our discussion regarding
possible extensions of Liouville's
theorem in the presence of negative curvature
begins with a result of Milnor \cite{Milnor}. A surface is said to be parabolic if 
it supports a 
nonconstant bounded harmonic function
and hyperbolic if it does not. So, a hyperbolic surface supports a bounded non constant harmonic
function, thus violating Liouville's theorem, whereas in a parabolic surface every bounded harmonic
function must be constant.
Milnor's result gives a curvature condition to distinguish between a parabolic and a 
hyperbolic surface. Next, we state a generalization of part of Milnor's result due to Greene and
Wu. In order to do so, we assume that $\mathbb{R}^n$ has a metric of the form
\begin{equation}
\label{eq:metric_general_form}
g=dr^2+\phi\left(r,\theta\right)d\theta^2,
\end{equation}
and we shall denote by $B_R\left(O\right)$ the open ball of radius $R>0$ centered at the origin.
In this case, recall that the curvature of $\left(\mathbb{R}^n, g\right)$ is given by
\[
K_g\left(r,\theta\right)=-\frac{\phi''\left(r,\theta\right)}{\phi\left(r,\theta\right)},
\]
where the prime ($'$) represents differentiation with respect to $r$.
\begin{theorem}
   Let $u$ be a harmonic function on a $\mathbb{R}^2$ with a metric of the form
   (\ref{eq:metric_general_form}), such that for $r\geq r_0>1$ its
   curvature function satisfies $K_g\left(r,\theta\right)\geq -\frac{1}{r^2\log r}$. 
   Let
   \[
   M\left(u;r\right)=\sup_{x\in B_R\left(O\right)}\left|u\left(x\right)\right|
   \]
   If 
   \[
   \liminf_{r\rightarrow \infty}\frac{M\left(u;r\right)}{\log\log r}=0,
   \]
   then $u$ is constant.
\end{theorem}
Therefore, a surface with a pole (that is, surfaces for which there
is a point $p$ where the exponential
map $\exp_p$ is a diffeomorphism), and whose curvature is $\geq -\frac{1}{r^2\log r}$ outside 
a compact subset, satisfies Liouville's theorem. A proof of this result can be found in \cite{Cortissoz2016}.

\medskip
Related to the previous result is the fact that in manifolds of negative curvature, Liouville's
theorem is flagrantly violated: not only do bounded nonconstant harmonic functions exists, but
even its asymptotic behaviour can be prescribed. In order to do this, a compactification
of $\mathbb{R}^n$, which we shall denote by
$\overline{\mathbb{R}^n}$ must be introduced (and which basically amounts
to adding the sphere at infinity), so that a Dirichlet problem at infinity
for the Laplacian can be posed as follows (see \cite{Choi}): Given $f\in C^0\left(\mathbb{S}^{n-1}\right)$,
find $u\in C^2\left(\mathbb{R}^n\right)\cap C\left(\overline{\mathbb{R}^n}\right)$
such that $\Delta u=0$ in $\mathbb{R}^n$ and $u=f$ on $\mathbb{S}^n$. 
Regarding the Dirichlet problem at infinity, in 1983, Anderson and Sullivan proved independently the following result.
\begin{theorem}[\cite{Anderson83, Sullivan83}]
    Let $M$ be a Cartan-Hadamard manifold whose sectional curvature is in between two negative constants,
    that is, there are $a,b\neq 0$ such that, the sectional curvatures $K$ of $M$ satisfy
    $-b^2\leq K\leq -a^2$. Then the Dirichlet problem is uniquely solvable for any 
    continuous $f$.
\end{theorem}

\medskip
In what follows, we shall give a generalisation of Liouville's Theorem on Riemannian cones (which
is part of an ongoing work with J. E. Bravo), and, in this regard, in the applications
we will see an interesting interaction between the warping function $\phi$ and the
first eigenvalue for the Laplacian $\Delta_{g_{\omega}}$ of $N$ when we are studying 
the Riemannian cones over $N$ (see below for the definition of a Riemannian cone).

\section{Preliminaries and Notation}

Here we review some of the theory and the results presented in \cite{Co2} that will be needed in the
following section.

\medskip
Given a closed $\left(n-1\right)$-dimensional Riemannian manifold $\left(N,g_{\omega}\right)$, we define the topological cone
\[
M=\left(\left[0,\infty\right)\times N\right)/\left(\left\{0\right\}\times N\right).
\]

Let $g$ be a metric defined on the cone $M$ by
\[
g=dr^2+\phi\left(r\right)^2g_{\omega},
\]
where $\phi\left(0\right)=0$ and $\phi'\left(0\right)=1$. We call $\phi$ the warping function,
and we shall call $\left(M,g\right)$
a Riemannian cone over $N$.
Observe that our definition of a Riemannian cone is a generalisation
of the concept of a metric cone, where $\phi\left(r\right)$ is taken as $r$.

\medskip
For $r>0$, the Laplace operator is defined as
\[
\Delta=\frac{\partial^2}{\partial r^2}+\left(n-1\right)\frac{\phi'}{\phi}\frac{\partial}{\partial r}+\frac{1}{\phi^2}\Delta_{g_{\omega}},
\]
where $\Delta_{\omega}$ is the Laplacian on $N$.
We say that $u\in C\left(M\right)\cap C^2\left(\left(0,\infty\right)\times N\right)$
is harmonic if $\Delta u=0$ on $\left(0,\infty\right)\times N$.

\medskip
Using separation of variables, we can produce harmonic functions in a cone as follows.
We denote by $f_{m,k}\left(\omega\right)$, $k=0,1,\dots,k_m$ the eigenfunctions
of the Laplace operator on $N$ for the eigenvalue $\lambda_m^2$, this is
\[
\Delta_{g_{\omega}} f_{m,k}=-\lambda_m^2 f_{m,k},\quad \lambda>0.
\]
We can thus write a harmonic function as 
\[
u_{m,k}=\varphi_m\left(r\right)f_{m,k}\left(\omega\right).
\]
In this case, it is easy to compute the equation that $\varphi_m$ must satisfy:
\[
\varphi_m''+\left(n-1\right)\frac{\phi'}{\phi}\varphi_m'-\frac{\lambda_m^2}{\phi^2}\varphi_m=0.
\]

\medskip
In what follows, we shall denote
\[
M_R=\left\{p=\left(r,\omega\right)\in M:\, 0\leq r\leq R\right\},
\]

Thus, given a continuous function
\[
h:N\longrightarrow \mathbb{R},
\]
it can be respresented a.e. by its Fourier expansion, say
\[
\sum_m\sum_k c_{m,k}f_{m,k}\left(\omega\right).
\]
Using this expansion, we can extend $h$ as a harmonic function to $M_R$ as
\[
u\left(r,\omega\right)=\sum_m\frac{\varphi_m\left(r\right)}{\varphi_m\left(R\right)}\sum_k c_{m,k}f_{m,k}\left(\omega\right).
\]
It can be shown that $u$ extends continously to $r=0$, and that this extension is unique.

\section{An abstract Liouville's theorem}

We come back now to Liouville's theorem. In order to state our next result, let us 
show a proof of Liouville's theorem devised by the author of this note and that
has been generalised by the author and by J. E. Bravo. Before we proceed, we want to
point out two facts: first, there is a simple proof of Liouville's theorem in the case
of an entire function using Fourier series; and, second, that in the case of dimension two,
the fact that harmonic functions remain harmonic after conformal deformation can be used to give another proof of Liouville's theorem as
presented below (see Section 3 in \cite{BravoCortissozPeters}): however, this 
second proof cannot be generalised to higher dimensions.

\medskip
So let's start with the promised proof.
Let $u:\mathbb{R}^2\longrightarrow \mathbb{R}$ be a harmonic function.
\[
M\left(u;R\right)=\sup_{x\in B_R\left(O\right)}\left|u\left(x\right)\right|
=\max_{x\in \partial B_R\left(O\right)}\left|u\left(x\right)\right|
\]
We let
\[
u_R\left(\omega\right)=u|_{\partial B_R\left(O\right)}.
\]
Since $u_R$ is smooth, we can expand it as the Fourier series
\[
u_R\left(\omega\right)=\sum_m \left(c_{m,R}e^{im\theta}+c_{m,R}^*e^{-mi\theta}\right),
\]
where the $*$ denotes complex conjugation.

The harmonic extension of $h_R$ to $B_R\left(O\right)$, and thus $u$, by uniqueness, is given by
\[
u\left(r,\omega\right)=\sum_m\left(\frac{\varphi_m\left(r\right)}{\varphi_m\left(R\right)}\right)
\left(c_{m,R}e^{im\theta}+c_{m,R}^*e^{-im\theta}\right),
\]
where
\[
\varphi_m\left(r\right)=\exp\left(\int_1^r\frac{m}{\phi\left(s\right)}\,ds\right).
\]
We will assume that $c_{0, R}=0$, which is equivalent to assume
that $u\left(O\right)=0$. Let us show that if $u\left(r,\theta\right)=o\left(\varphi_1\left(r\right)\right)$ then
$u$ must be constant.

Fix $R_0\ll R$. Given $\epsilon>0$,
we let $L$ be such that 
\[
\left\|\sum_{\left|m\right|\geq L}\left(\frac{\varphi_m\left(r\right)}{\varphi_m\left(R\right)}\right) 
\left(c_{m,R}e^{im\theta}+c_{m,R}^*e^{-im\theta}\right)\right\|^2_{L^2\left(B_{R_0}\left(O\right)\right)}\leq \epsilon.
\]

Next, we estimate $\left\|u\right\|_{L^2\left(B_{R_0}\left(O\right)\right)}$:
\begin{eqnarray*}
\left\|u\right\|_{L^2\left(B_{R_0}\left(O\right)\right)}^2&\leq&2\int_0^{R_0} \int_{\mathbb{S}^1}
\left|\sum_{\left|m\right|< L} c_{m,R}\left(\frac{\varphi_m\left(r\right)}{\varphi_m\left(R\right)}\right) e^{im\theta}\right|^2
\phi\left(r\right)\,
d\omega\,dr+\epsilon\\
&=&
2\int_0^{R_0}\sum_{m<L}\left(\frac{\varphi_m\left(r\right)}{\varphi_m\left(R\right)}\right)^{2}
\left|c_{m,R}\right|^2\,\phi\left(r\right)\,dr+\epsilon.
\end{eqnarray*}
Therefore, if say $R\geq 2R_0$, we can estimate
\begin{eqnarray*}
\left\|u\right\|_{L^2\left(B_{R_0}\left(O\right)\right)}^2&\leq&
\frac{1}{\varphi_1\left(R\right)^2}\int_0^{R_0}\sum_{m<L}\left(\frac{\varphi_m\left(R_0\right)}{\varphi_m\left(R\right)}\right)^2
\varphi_1\left(R\right)^2
\left|c_{m,R}\right|^2\,r\,dr+\epsilon\\
&\leq&
\frac{C}{\varphi_1\left(R\right)^2}\int_0^{R_0} \sum_{m<L} \left|c_{m,R}\right|^2\,\phi\left(r\right)\,dr+\epsilon,
\end{eqnarray*}
using the fact that $h_R=o\left(\varphi_1\right)$, we then have
\begin{eqnarray*}
&\leq&\frac{C}{\varphi_1\left(R\right)^2}o\left(\varphi_1\left(R\right)^2\right)\int_0^{R_0}\,\phi\left(r\right)\,dr+\epsilon\\
&\leq& \frac{C_{R_0}}{\varphi_1\left(R\right)^2}o\left(\varphi_1\left(R\right)^2\right)+\epsilon,
\end{eqnarray*}
and hence, if $R>0$ is large enough,
\[
\left\|u\right\|_{L^2\left(B_{R_0}\left(O\right)\right)}^2 < 2\epsilon.
\]
Being $\epsilon>0$ arbitrary, we must have
\[
\left\|u\right\|_{L^2\left(B_{R_0}\left(O\right)\right)}^2=0,
\]
and thus $u\equiv 0$.

\medskip
As mentioned before, the previous argument inspired the following result, whose proof follows along the same lines as shown above, and which shall 
be presented elsewhere \cite{BravoCortissoz}.
\begin{theorem}
\label{thm:main_theorem}
    Assume there is an $A\left(R\right)\rightarrow \infty$ such that $\varphi_m\left(R\right)\gg A\left(R\right)$,
    $m\geq 1$. Given a harmonic
    function $u:M\longrightarrow \mathbb{R}$, let
    \[
    u_R\left(\omega\right)=u\left(R,\omega\right).
    \]
    If
    $u_R\left(\omega\right)=o\left(A\left(R\right)\right)$, then $u$ is constant.
\end{theorem}

\section{Applications}

Let us show a couple of applications of Theorem \ref{thm:main_theorem}. Although the results we shall
show are not sharp, they demonstrate the applicability of our methods.

\medskip
First, let us assume that the warping funtion satisfies $\phi'\geq 1$ (which occurs whenever the
radial curvature is nonpositive). Then, the author and J. E. Bravo have shown that for $n\geq 3$ we can use as
$A\left(r\right)$ the following function (this result is not sharp)
\[
A\left(r\right)=\exp\left(\lambda_1^2\min\left\{\frac{1}{\sqrt{2}\lambda_1},\frac{1}{2}\right\}\int_1^r \frac{1}{\phi^{n-1}}\left(\sigma\right)
\int_{0}^{\sigma}\phi^{n-3}\left(\tau\right)\,d\tau\,\, d\sigma\right),
\]
where $\lambda_1^2$ is the first nontrivial eigenvalue of the Laplace operator defined in $\left(N,g_{\omega}\right)$,
that is, there is a nonzero smooth function $v$ such that
\[
-\Delta_{g_{\omega}} v=\lambda_1^2 v, \quad \lambda_1>0.
\]
Therefore, if 
\[
\int_1^r \frac{1}{\phi^{n-1}}\left(\sigma\right)
\int_{0}^{\sigma}\phi^{n-3}\left(\tau\right)\,d\tau\,\, d\sigma\rightarrow\infty
\quad \mbox{as}\quad r\rightarrow \infty,
\]
the cone does not support a bounded nonconstant harmonic 
function. Let us work out the case of $\phi\left(r\right)=r$ and $N=\mathbb{S}^{n-1}$
with an arbitrary metric. In this case
\begin{eqnarray*}
    A\left(r\right)&=&\exp\left(\frac{1}{n-2}\lambda_1^2\min\left\{\frac{1}{\sqrt{2}\lambda_1},\frac{1}{2}\right\}\ln r\right),
\end{eqnarray*}
and hence, if $\lambda_1\geq \sqrt{2}$,
we obtain:
\[
A\left(r\right)=r^{\frac{\lambda_1}{\sqrt{2}\left(n-2\right)}},
\]
and this gives a bound from below for the slowest possible growth of a nonconstant  harmonic function defined 
on the cone with metric $dr^2+\phi\left(r\right)^2g_{\omega}$, $g_{\omega}$ being an arbitrary metric
on $\mathbb{S}^{n-1}$. If $g_{\omega}$ is the round metric, then it is well known that $\lambda_1=\sqrt{n-1}$, and
thus a bound from below for the slowest possible growth of
a nonconstant harmonic function in this case would be, for $n\geq 3$,
$r^{\frac{\sqrt{n-1}}{\sqrt{2}\left(n-2\right)}}$, which does give Liouville's theorem in its classical version,
but not its sublinear version (with the remarkable exception of $n=3$). 

\medskip
On the other hand, if the radial curvature is nonnegative, then, the expression above can be simplified a bit 
and we can take (see \cite{BravoCortissoz})
\[
A\left(r\right)=\exp\left(\lambda_1^2\min\left\{\frac{1}{\sqrt{2}\lambda_1},\frac{1}{2}\right\}\int_1^r\frac{1}{\phi\left(s\right)}\,ds\right).
\]
As a consequence, if $\phi\left(r\right)=o\left(r\right)$ then, on a cone with nonnegative radial curvature,
there are no nontrivial harmonic functions of polynomial growth. This should be compared with the
example of the paraboloid obtained by rotating the curve $y=x^2$ around the $y$-axis (see \cite{BravoCortissozPeters}),
for which the result just described holds.

\end{document}